\documentclass[preprint,12pt]{article}

\usepackage{setspace}
\usepackage{pdfpages}

\usepackage{amsthm,amsmath,amssymb,amsfonts,amscd,amsbsy}
\usepackage{graphicx}
\usepackage[ruled,noline,linesnumbered]{algorithm2e}
\usepackage[tight]{subfigure}
\usepackage{fullpage}
\usepackage{color}
\usepackage{url}
\usepackage{mdframed}

\DeclareGraphicsExtensions{.jpg,.png,.pdf}
\graphicspath{{./fig/}}

\def\calA{\mathcal{A}}

\def\calD{\mathcal{D}}

\def\calF{\mathcal{F}}

\def\calM{\mathcal{M}}

\def\calP{\mathcal{P}}

\def\calS{\mathcal{S}}

\def\RR{\mathbb{R}}
\def\SS{\mathbb{S}}

\newcommand{\ttt}[1]{\texttt{#1}}

\newcommand{\bmat}[1]{\left(\begin{array}{#1}}
\newcommand{\emat}{\end{array}\right)}

\title{Project Artie: An Artificial Student for Disciplines Informed by Partial Differential Equations\thanks{Submitted to Computer Applications in Engineering Education}}
\date \today
\author{Anthony T Patera \thanks{Department of Mechanical Engineering, Massachusetts Institute of Technology \newline \mbox{}\hspace{.21in} Room 3-266, MIT, 77 Massachusetts Avenue, Cambridge, MA 02139 USA \newline \mbox{}\hspace{.20in} patera@mit.edu} }

\begin{document}

\maketitle

\begin{abstract}

We present an Artificial Student, ``Artie,'' for engineering science disciplines in which the mathematical model is a partial differential equation (PDE); Artie considers here the particular case of steady heat conduction. Artie accepts problem statements posed in natural language. Artie provides a symbolic-numeric approximate solution: the PDE field; scalar Quantities of Interest (QoI), expressed  as functionals of the field. The problem statement will typically not provide explicit guidance as to the equation or approximations which should be invoked. We also present Artie+, who provides the finite element solution to the PDE: the exact solution to within a prescribed  tolerance controlled by an {\it a posteriori} error estimator.

Artie comprises four technical ingredients: Natural Language Processing: We proceed in two stages:  domain-independent Google Natural Language syntax analyzer followed by frame-specific conduction parser. PDE Template: The PDE is exploited by the conduction parser  to extract geometry, boundary conditions, and coefficients; subsequent approximations are deduced from this ground-truth description. Problem Classes, Geometry Classes; Components, Systems: A problem class places requirements on spatial domain, boundary conditions, properties, and QoI; associated to each problem class are several geometry classes. A component is an instantiation of the geometry class for prescribed  geometric and PDE parameters; a system is represented as an assembly of connected components. Variational Formulation: We consider the weak statement and minimization principle to formulate the PDE and develop suitable approximations; implementation proceeds through static condensation and direct stiffness assembly over component ports.

We describe and illustrate a prototype implementation of Artie and Artie+.
\end{abstract}

\section{Motivation}

We present the software in the guise of an artificial student, and it is thus appropriate to briefly summarize our perspective on undergraduate education  in the engineering sciences.

\subsection{Pedagogical Background}
\label{subsec:PB}

Traditional fields of continuum mechanics, such as conduction heat transfer or linear elasticity, are well understood: constitutive laws and conservation principles yield a well-posed governing partial differential equation (PDE). Nevertheless, education in these  disciplines remains challenging: the abundance of behaviors for different geometries, materials, and boundary conditions is matched only by the paucity of corresponding closed-form (hence indisputable) solutions.

Engineers and engineering students hence typically --- or at least classically --- replace the governing PDE with a look-up table: the dictionary comprises a few canonical problem classes and associated closed-form approximation procedures; any given problem is modified such that the perturbed problem conforms to a ``nearest'' canonical problem class; the perturbed problem is then amenable to closed-form approximate solution. We can view the process conceptually as a Voronoi tesselation of problem space in which the canonical problem classes serve as generators. The resulting (symbolic closed-form) approximate problem solutions are very useful for prediction and design: transparent, accessible, and rapidly evaluated.

Numerical solution of the governing PDE of course offers much higher accuracy in particular for problems far from any canonical problem class. Not suprisingly, PDE computation now plays the central role in industry research, development, and design. Yet we continue to teach (implicitly) the Voronoi approach in most of our undergraduate engineering subjects, and for good reason: re-purposed closed-form approximations retain relevance even in the digital era, most notably for (i) conceptual and preliminary design --- motivated by the expense and relatively slow response of PDE solvers, and  (ii) verification of numerical computations --- motivated by the opacity of PDE solvers and associated software. Closed-form approximations are particularly valuable for the detection of blunders committed in the preprocessing, processing, or postprocessing stages of simulation.

\subsection{The Artificial Student} 

We introduce in this paper an Artificial Student, ``Artie,'' for engineering science disciplines in which the governing mathematical model is a partial differential equation; in this first embodiment, Artie considers the particular case of undergraduate-level steady heat conduction. 
\begin{itemize}
\item[] Artie accepts problem statements as posed in natural language to our (actual) students. {\em In this first paper, we consider natural language processing but not image processing; the problem statement text is suitably elaborated to provide geometric information.} 
\item[] Artie provides an approximate problem solution: (i)  the PDE field --- here, the temperature; (ii) scalar Quantities of Interest (QoI), expressed as functionals of the PDE field  --- here (say), the heat transfer rate over a prescribed surface. {\em In this first paper we shall only consider problem statements which conform to some problem class, and hence we eliminate the source of error associated with problem perturbation.}
 \end{itemize}
Note we consider only analysis (or forward) problems, not inverse or design problems.
 
We emphasize that the problem statement will typically not provide explicit guidance as to the particular equation or approximations which should be invoked: the student, and Artie, must deduce an appropriate approximation procedure based on the stated geometry, boundary conditions, physical properties, and QoI. We also note that the actual student and Artie need not, and typically will not, pursue the same inference or approximation procedure; however, we do require that the actual student and Artie will each arrive ultimately at a valid solution --- and typically the same solution --- to any given problem within a particular problem class.

At present Artie is a rather rudimentary prototype. But ultimately, upon subsequent development and elaboration, we envision several roles for Artie:

\begin{enumerate}
\item (education) Artie can provoke. A successful Artie, subject to the same grading standards as our actual students, should receive an ``A'' in an undergraduate engineering science subject. Should we thus reconsider what we teach our students? Or how we assess our students?

\item (education) Artie can provide learning assistance. Artie can reveal to students intermediate results in symbolic form. The student can thus identify not only a flaw in their approach but also possibly the source of the error. (Artie can not necessarily provide the student with {\em procedural guidance}, since the actual student and Artie may pursue different albeit equivalent approximation approaches.)  

\item (education) Artie can provide teaching assistance. Artie can serve as ``surrogate student'' in the development of problems for assignments and exams. In particular, Artie can detect problem statements which are ambiguous or misleading; Artie can also identify (underlying) mathematical models which are ill-posed.

\item (professional practice) Artie can provide blunder detection. Artie's closed-form approximate solutions can serve to identify blunders in PDE numerical solution procedures: incorrect input specifications, insufficiently refined discretizations, and errors in interrogation. (We emphasize ``blunders'':  the error in Artie's approximate solution must be sufficiently less than the error in the numerical solution.)

\end{enumerate}
Artie can also serve, through the development process, in the abstraction of approximation procedures.

In this paper we also describe Artie+. Given a problem statement, Artie+ provides Artie's closed-form symbolic-numeric approximate problem solution but also the finite element (FE) solution to the PDE. The latter is the exact solution to the problem statement within a prescribed tolerance as enforced by an asymptotic {\it a posteriori} error estimator: we control the error in the PDE field in an appropriate energy norm; we also directly control  the error in the QoI. In this paper we impose tight tolerances such that the FE solution is effectively exact.

Our prototype of Artie+ is even more rudimentary than our protype of Artie, and in particular we address at present only problems defined over two-dimensional domains (hence infinite or insulated in the third direction). However, ultimately, we envision several roles for Artie+:

\begin{enumerate}
\item (education) Artie+ can provide supplemental  information and perspective. Artie+ can confirm --- or refute --- the relevance of Artie's closed-form approximate solution through direct comparison with the exact (PDE numerical) solution. Artie+'s rendering of the PDE field can furthermore highlight the source of any substantial discrepancies.

\item (professional practice) Artie+ can provide more flexible and more rigorously certified engineering analyses. Artie's closed-form approximations can serve for rapid evaluation; Artie+'s PDE numerical solution can be invoked, as necessary, for confirmation. (Note that since the closed-form prediction and the PDE numerical solution share data and code, Artie+ is not necessarily reliable as a blunder detector.)
\end{enumerate}
We also note that, in some sense, Artie+ provides a natural language pre-processor for PDE numerical solution.
\section{Technical Approach}

\subsection{Principal Ingredients and Related Work}

There are four main ingredients to Artie and Artie+:

\begin{enumerate}
\item {\em Natural Language Processing}. We consider a two-step approach to syntax analysis similar to the architecture described in \cite{twopart_parser}: in the first stage we apply the general (domain-independent) Google Natural Language syntax analyzer \cite{Google_NLP} to identify tokens, parts-of-speech, and verb tense; in the second stage we apply a domain-specific syntax analyzer of our own conception, which we denote a (steady) {\em conduction parser}. The latter treats any given problem statement in terms of a heat conduction ``frame''  of reference \cite{frames1,frames}. 

At present, Artie has no real intelligence: the conduction parser learns from instances only through the intermediary of the code developer, who expands and modifies the parsing algorithm as needed for each new problem statement encountered; in principle, increasingly few (and ultimately no) modifications are required as the training suite is expanded. In actual fact, our interest is less in natural language than in non-unnatural language: Artie can still well serve the purposes cited above even if we restrict syntax, and in particular discourage complicated sentence structures apt to confuse; the latter are arguably poor scientific and pedagogical practice in any event.

\item {\em PDE Template}. The conduction parser is certainly informed by the particular context and language of conduction heat transfer \cite{Lienhard}. But more generally the conduction parser is informed by the underlying PDE which describes heat conduction: we take the viewpoint that any rational approximation must ultimately derive from the ground truth; we thus (implicitly) complete a PDE template as the first and crucial step in the approximate problem solution process. The PDE Template corresponds to data structures which characterize the problem spatial domain $\Omega$, the type of boundary conditions over different parts of the spatial domain boundary $\partial \Omega$, the surface source terms, any volumetric source terms, the coefficients which appear in the weak statement of the PDE, and the type of QoI; these specifications are collectively constrained to ensure a well-posed problem \cite{QandV}. We note that in most cases the template considers the abstract mathematical form such that Artie would require relatively little modification to treat different physical disciplines governed by (currently) second-order elliptic PDEs. For example, the PDE Template boundary condition types Dirichlet, Neumann, and Robin correspond respectively to the problem statement specifications of ``temperature,'' ``flux,'' and ``heat transfer coefficient.''

In this context we mention two earlier related efforts. We first discuss the FEniCS project \cite{FEniCS}. The goals of Artie+ are similar to the goals of FEniCS: a natural language interface for description and ultimately solution of PDEs. However, Artie considers the natural language of engineers, a modestly expanded form of English, whereas FEniCS considers the natural language of mathematicians --- a precise mapping from symbols to interpretation. Conversely, in \cite{alg_wordproblems}, the authors do consider solution of mathematical problems posed in (truly) natural language: \cite{alg_wordproblems} addresses word problems which yield small systems of algebraic equations. In fact, perhaps surprisingly, Artie's task is simpler than the goals pursued in \cite{alg_wordproblems}: the identity of our variables is readily extracted from the specific (heat conduction) context.

\item {\em Problem Classes and Geometry Classes; Components and Systems}. As described earlier, Artie and Artie+ will only accept a problem statement which belongs to an anticipated {\em problem class}: the criteria for membership involve the geometry and topology, the boundary conditions, the physical properties, and the QoI. Associated to each problem class we define several parametrized {\em geometry classes}. 

The problem class and geometry class are coupled through the notion of a {\em component}: a component is an instantiation of the geometry class for prescribed values of the geometric parameters and the (local restriction of the) PDE parameters. Components provide several advantages: components are a convenient fashion by which to construct systems with discontinuous coefficients, as arise frequently in engineering analysis; components often refer to actual building blocks encountered in engineering practice. Components are endowed with (i) faces for assignment of boundary conditions and connection with other components, and (ii) ports for representation of degrees of freedom. A {\em system} associated to a particular problem statement is an assembly of components connected at compatible faces. 

At present, Artie treats two problem classes, which indeed are sufficient to treat a majority of undergraduate steady heat conduction questions:

\begin{itemize}
\item[I] Quasi-1d Systems.  This {\em problem class} places requirements on geometry and parameters: geometry in (the axial coordinate) $x$ is homogeneous such that the cross-sectional area $A$ and perimeter $P$ are constant; the Biot number, Bi $\equiv h_{\max} (A/P)/k_{\min}$, is small compared to unity, where $h_{\max}$ and $k_{\min}$ correspond to the maximum heat transfer coefficient and minimum thermal conductivity of the system; all standard boundary conditions are supported; the QoI may take the form of temperature anywhere in the domain, and flux or heat transfer on either axial face of $\partial \Omega$.  We associate to this problem class a single {\em geometry class}: right cylinder with specified cross-section shape; at present, we consider only rectangular cross section, with two geometric parameters, however more cases can readily be included by minor expansion of the class definition. A {\em component} has (i) three faces, two rectangular axial faces and a lateral face, and (ii) two ports, which correspond to the axial faces. This problem class can be readily extended to geometry which is no longer homogeneous in $x$ but only slowly varying in $x$, however we must then relax ``closed-form'' to include numerical solution of {\em ordinary} differential equations.

Most notably, quasi-1d systems are relevant to one-dimensional walls and, more importantly, thermal fins \cite{Lienhard}; the latter constitute a highly relevant example of extended-surface heat transfer which is furthermore within undergraduate reach as regards analysis and application.

\item[II] Generalized Walls. This {\em problem class} places requirements on geometry, boundary conditions, and QoI: the spatial domain $\Omega$ must be represented as a conforming union of bricks; all exposed surfaces are insulated except for the faces at $x_1 = x_1^{\rm left}$ and $x_1 = x_1^{\rm right}$ exposed to respective (heat transfer coefficient, fluid temperature) pairs ($h^{\rm left}$, $T^{\rm left}$) and ($h^{\rm right}$, $T^{\rm right}$); the QoI is a nondimensional heat transfer rate over all  faces at $x_1 = x_1^{\rm left}$. Here $x_1^{\rm left} = \min_{x_1 \in \bar{\Omega}}$ and $x_1^{\rm right} = \max_{x_1 \in \bar \Omega}$ define the smallest slab (in $x_1$) which contains  $\Omega$. We associate to this problem class a single {\em geometry class}: a parallelpiped with specified dimensions in $x_1$, $x_2$, and $x_3$. A {\em component} has (i) six rectangular faces, and (ii) two rectangular ports, which correspond to the $x_1 = \text{ const}$ faces, which we will denote simply ``$x_1$ faces''; we introduce only two ports and hence two degree of freedom per component in anticipation of the subsequent closed-form (1d) approximation. We note that this problem class can be readily extended to consider appropriate combinations of temperature and flux boundary conditions at $x_1 = x_1^{\rm left}$ and $x_1 = x_1^{\rm right}$.

Most notably, generalized walls are relevant to actual walls in particular intended to insulate. The generalization to brick construction then permits heterogeneous materials and composites, voids, and more complex geometries, and furthermore illustrates the effects of geometric contraction and expansion on heat flow.

\end{itemize}
We can certainly view Artie as an expert system, or more precisely a collection of expert systems associated with respective problem classes. Artie differs from current heat transfer expert systems (typically focused on heat exchangers, for example \cite{Afgan,Cochran}) in a number of ways: Artie admits inputs in natural language; Artie's problem classes are defined by a very high-dimensional parameter space; Artie's set of ``rules'' is  relatively complicated and often implicit; and finally, Artie relies relatively weakly on a knowledge or experiential database (we discuss the latter in Section \ref{subsec:CP}). But we must also emphasize that Artie treats academic undergraduate problems which are very simple compared to the industrial heat exchanger problems considered in \cite{Cochran,Afgan}.

\item {\em Variational Formulation and Approximation}. We consider the weak statement and associate minimization principles both to formulate the PDE and also to develop suitable approximations \cite{QandV}.  The variational formulation offers several advantages: natural boundary conditions are very easily accommodated; systematic approximation is possible by consider of appropriate (sub)spaces; the minimization principle offers in many cases useful lower and upper bound constructions \cite{bounds_us,bounds_Elrod}, as discussed further in Section \ref{sec:Ex} example \texttt{wall-3d}; and finally, we can pursue direct stiffness assembly \cite{Bab} at the component level to automatically form the system matrix. We note that prior to direct stiffness assembly we perform (symbolically) static condensation \cite{Wi_IJNME} to restrict the degrees of freedom to the ports.

As indicated earlier, actual students and Artie will typically follow different procedures to derive the temperature field and associated QoIs. In particular, whereas actual students might invoke the strong form, thermal resistance concepts, conservation of energy (network Kirchoff Laws), and equivalent resistance, Artie pursues a variational formulation with suitable approximation spaces. It would be possible to enhance Artie to translate variational results (after the fact) into ``student form.'' 
\end{enumerate}

\subsection{Conduction Parser}
\label{subsec:CP}

We provide a few details of the conduction (PDE) parser, in particular to demonstrate the relatively generic nature of the procedure. We recall that the conduction parser takes as input the Google Natural Language syntax analysis: the tokens, parts-of-speech, and verb tenses associated with the problem statement. In actual practice, we call the Google Natural Language parser several times in order to disambiguate the original problem statement.

\begin{enumerate}
\item {\em Syntax Preparation}. The conduction parser first modifies the Google syntax analysis to reflect the technical nature of our frame. As but one example, ``normal'' can be a noun, not just an adjective, in engineering science analysis. We also introduce our own escape characters to signal symbol and equation delimiters and also symbol and equation tokens.

\item {\em Entities, Snippets, and Attributes}. We next identify {\em entities}, compound nouns (typically) followed by adpositional phrases. We then isolate {\em snippets}: simple subject-predicate-object (or -complement) sequences expressed in terms of entities and a present-tense verb. We then derive, from the snippets, {\em attributes} associated with each entity. We also perform some disambiguation operations, in particular given that we must perform the syntax analysis for the entire problem statement and not just individual sentences;  for example, ``handle of the spoon" and ``spoon handle'' refer to the same entity and are coalesced. 

\item {\em Commonsense Incorporation}. We then perform a similar syntax analysis not on the problem statement but rather on a ``commonsense database'' which serves all problems in all problem classes.  The latter includes information, largely non-technical, which we would expect an actual student to know, for example ``air is a gas", or ``a spoon is a solid object". We then append commonsense (entity, entity attributes) as problem statement entity attributes according to certain subset and  inclusion rules; for example, commonsense sentence ``A spoon is a solid object.'' is incorporated as problem statement (entity = ``spoon'', attribute = ``a solid object''); the commonsense syntax tree is then discarded. 

\item {\em Solid and Fluid Entities}. We can now identify solid entities and fluid entities, a subset of which will ultimately define our components and then degrees of freedom. In fact, Artie can also deduce solid or fluid ``state'' from other diagnostics, such as connection; in future we will take advantage of multiple predictors to provide more robust inference.

\item {\em Inheritance and Instantiation}. We next look for ``inheritance'' words to identify parent-child relationships; for example, ``a spoon consists of a head and handle'' yields inheritance pairs ``spoon-handle'' and ``spoon-head''. These inheritance pairs then serve to pass attributes from parent to children; for example, ``head'' and ``handle'' are now labeled as solid objects even though ``head'' in isolation is not evidently a solid object. 

In a similar fashion we search for ``instantiation'' indicators, for example ``each'', to identify archetypes. We then perform a second attribute incorporation procedure now over the problem statement entities; ``spoon'' and ``spoon geometry'' are separate entities, however we associate ``spoon geometry'' attributes to ``spoon'' as well; even more importantly, attributes of (say) archetype ``brick'' --- for example, geometry class as identified by particular keywords --- are then automatically adopted by instantiations ``brick 1'', ``brick 2'', ``brick 3''. 

In engineering analysis we often construct systems through assembly and instantiation, and it would be laborious or even impossible --- and distinctly un-natural --- to  individually specify attributes clearly shared by many parts.

\item{\em Connection}. We next search for ``connection'' words, for example ``a wall {\em separates} inside air and outside air'', to construct a system graph: the nodes are entities, and the edges constitute possible heat transfer paths. In future the system graph might be derived from a figure and image processing, however at present we rely on text analysis. 

\item {\em Component Identification}. A  component is defined as an entity which satisfies the following four conditions: (i) the entity --- interpreted as a node of our connection graph --- is of degree $> 0$, (ii) the entity is solid, (iii) the entity is {\em not} an ``insulator'' or made of insulating material, and (iv) the entity is {\em neither} an inheritance parent entity {\rm nor} an instantiation archetype entity.

\item {\em Coordinate Variables and Spatial Domains}. We next extract coordinate variables and the spatial domains associated with each component; we can then construct the system spatial domain through component face connectivity information. Note all processing is symbolic, and hence there are no issues related to numerical precision in the identification of topology and geometry. We can subsequently develop the usual local-to-global mapping from component and local port to unique global degree of freedom (associated to coincident local ports). In future, spatial domains may be deduced or corroborated from a figure, but at present we are limited to text analysis.

\item {\em Boundary Conditions and QoI}. We then identify boundary conditions  through structured searches of entity attributes and snippets informed by appropriate keywords --- ``temperature'', ``heat flux'', and ``heat transfer coefficient''. Note a boundary condition such as heat transfer coefficient requires not just a heat transfer coefficient but also a fluid temperature, each of which might appear in different sentences within the problem statement.  QoI, such as temperature, heat flux, and heat transfer rate, may be found in a similar fashion with additional criteria related to ``find'' words.

\item {\em Physical Properties}. For steady conduction the only relevant physical property is thermal conductivity. We can identify thermal conductivity either directly (as in the examples in this paper) or through attributes which contain ``material'' information; in future, we will also include a material properties database in the same fashion as the commonsense database.

\item {\em Parameter Values}. Finally, we search for equalities which do not involve the coordinate variables in order to isolate the numerical parameter values associated to our symbolic variables; the former are important (i) to determine parameter regime and hence problem class, quite independent of numeric evaluation of the solution, but also (ii) to instantiate numerical solutions for particular physical systems of interest. 

\end{enumerate}
At this stage all information is available for approximate solution of the PDE (Artie) and FE solution of the PDE (Artie+). Note the approximate solution is expressed in terms of a (static-condensation) basis directly deduced from the problem class and component specifications; the corresponding stiffness matrix and load vector are then easily evaluated symbolically, and the system matrix formed from the local-to-global mapping.

We conclude with several remarks. First, with only few and small exceptions, most of the preceding steps are not overly specific to heat conduction, and could be readily extended to other continuum disciplines. Second, as already indicated, the connection graph and component spatial domains would typically be deduced by (actual) students from a figure. Absent image processing at present, we resort to a verbal explanation. In fact, Artie prepares from the verbal explanation figures with the connection graph and geometry, and thus we could reserve the graph and geometry text solely for Artie as pre-processor and then supply the students with the corresponding figures. Ultimately, however, Artie must graduate to include image processing and character recognition as well. Third, and finally, in this first embodiment of Artie we do not explicitly include units: all quantities are given implicitly in SI units (note temperature is in degrees Celsius, not Kelvin). In the future, we will include explicit units in the problem statements, however in this first embodiment we prefer to test Artie without the rather heavy-handed clues provided by units; Artie's performance will only improve once we can include units in our inference procedure. In fact, units are often only associated with numerical values, and hence units should only serve as corroboration.

\section{Examples}
\label{sec:Ex}

We precede the examples with a few details related to implementation. Artie and Artie+ are implemented in MATLAB \cite{matlab}: MATLAB supports the string, symbolic, and numeric processing required by Artie and Artie+. A Google Natural Language client is hosted within MATLAB \cite{urlread2,PhGNL} to permit ready and interactive query. In principle Artie can provide symbolic results for the PDE field and QoI; in practice, the MATLAB-produced symbolic expressions are not tractable or enlightening for more than a single component, and thus we typically resort to numerical evaluation of our symbolic component stiffness matrices prior to assembly and inversion. The finite element (FE) procedures of Artie+ are implemented within MATLAB by the \texttt{fem2d} template package \cite{Masafem2d} which in turn relies on the \texttt{DistMesh} \cite{Per} mesh generator; we consider  simple two-level (asymptotic)  error estimators \cite{Bab} which inform a \texttt{fem2d} Nearest Vertex Bisection (NVB) adaptive refinement procedure.\footnote{The Artie and Artie+ software is available ``as-is'' free-of-charge under an open-source license. All required third-party packages except the Google Natural Language Processor are also available free-of-charge under an open-source license; the Google Natural Language Processor, for Artie's modest requirements, is very inexpensive. We emphasize that Artie is a research code not intended (in the present version)  for production within the context of actual subject deployment.}

We now turn to three examples. The first two examples exercise problem class Quasi-1d; the third example illustrates problem class Generalized Wall. In all cases we consider either heat transfer coefficient or zero-flux boundary conditions, though Artie is also instrumented to treat temperature and non-zero flux boundary conditions.  We present here the problem statements in latex format, but for full disclosure we include in the appendix the plain-text (verbatim) input files actually processed by Artie; the latter are adorned with our particular escape characters and any special allowances related to the Google Natural Language Processor or subsequent MATLAB symbolic manipulation.


\subsection{Example 1: \texttt{wall-1d}}. 

We present the problem statement:

\begin{itemize}
\item[]
A composite wall separates inside air and outside air. The inside air is maintained at temperature $T_{\rm in}$; the outside air is maintained at temperature $T_{\rm out}$. The composite wall comprises three layers: a fir layer, a pine layer, and a cedar layer; the pine layer is in between the fir layer and the cedar layer. The fir layer is exposed to the inside air; the cedar layer is exposed to the outside air. Figure 1 is Artie's depiction of the heat transfer paths.

The composite wall is a right cylinder with rectangular cross-section of dimensions $a$ and $b$; the coordinate through the wall is $x$. The fir layer is of length $L_{\rm f}$; the pine layer is of length $L_{\rm p}$; the cedar layer is of length $L_{\rm c}$. The spatial domain of the fir layer is $0 < x < L_{\rm f}$; the spatial domain of the pine layer is $L_{\rm f} < x < L_{\rm f} + L_{\rm p}$; the spatial domain of the cedar layer is $L_{\rm f} + L_{\rm p} < x < L_{\rm f} + L_{\rm p} + L_{\rm c}$. Figure 2 is Artie's construction of the wall geometry.

Let $h_{\rm in}$ denote the heat transfer coefficient from inside air to fir layer prescribed over the  face at $ x = 0 $; let $h_{\rm out}$ denote the heat transfer coefficient from cedar layer to outside air prescribed over the face $ x = L_{\rm f} + L_{\rm p} + L_{\rm c} $. The fir layer, pine layer, and cedar layer are insulated on the lateral faces.

The thermal conductivity of the fir layer is $k_{\rm f}$; the thermal conductivity of the pine layer is $k_{\rm p}$; the thermal conductivity of the cedar layer is $k_{\rm c}$.

Plot the temperature distribution as a function of $x$. You may use the following parameter values: $T_{\rm in} = 23$, $T_{\rm out} = 0$, $a = 0.1$, $b = 0.1$, $L_{\rm f} = 0.05$, $L_{\rm p} = 0.1$, $L_{\rm c} = 0.05$, $h_{\rm in} = 10$, $h_{\rm out} = 100$, $k_{\rm f} = 0.2$, $k_{\rm p}= 0.1$, and $k_{\rm c} = 0.05$. $\Box$
\end{itemize}
We observe that this problem is in problem class Quasi-1d, and in fact the problem as stated is exactly one-dimensional. The geometry class is right cylinder with rectangular cross-section; in practice, we can provide aliases and hence refer to  ``slab'' rather than right cylinder. (The geometry description is mathematically informal, but context (problem class) serves to fill in the blanks.) We note that ``plot'' is mandatory within problem class Quasi-1d, and hence the first sentence of the last paragraph is parsed by Artie but no inference is required. We do not request a QoI.

\begin{figure}[h!]
\begin{center}
 \includegraphics[scale=0.24]{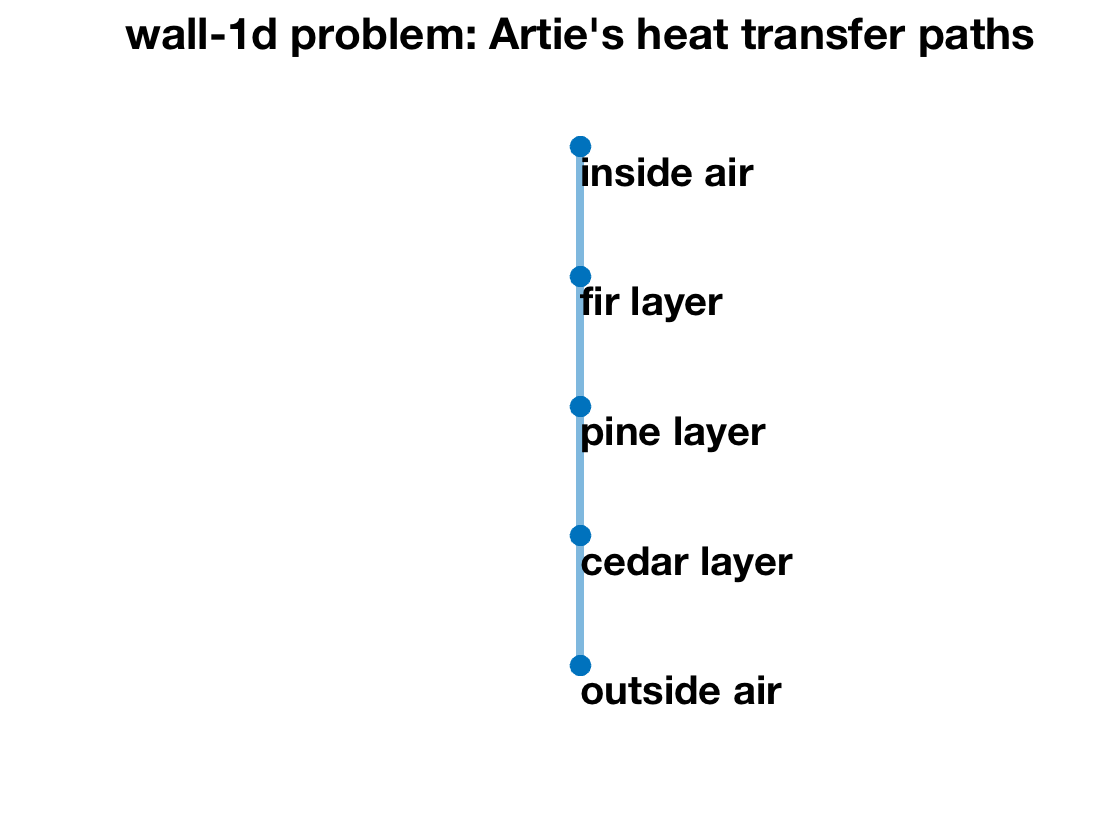}
\end{center}
\caption{Artie's depiction of the heat transfer paths for problem \texttt{wall-1d}.}
\label{fig:wall-1d_fig1}
\end{figure}

\begin{figure}[h!]
\begin{center}
 \includegraphics[scale=0.3]{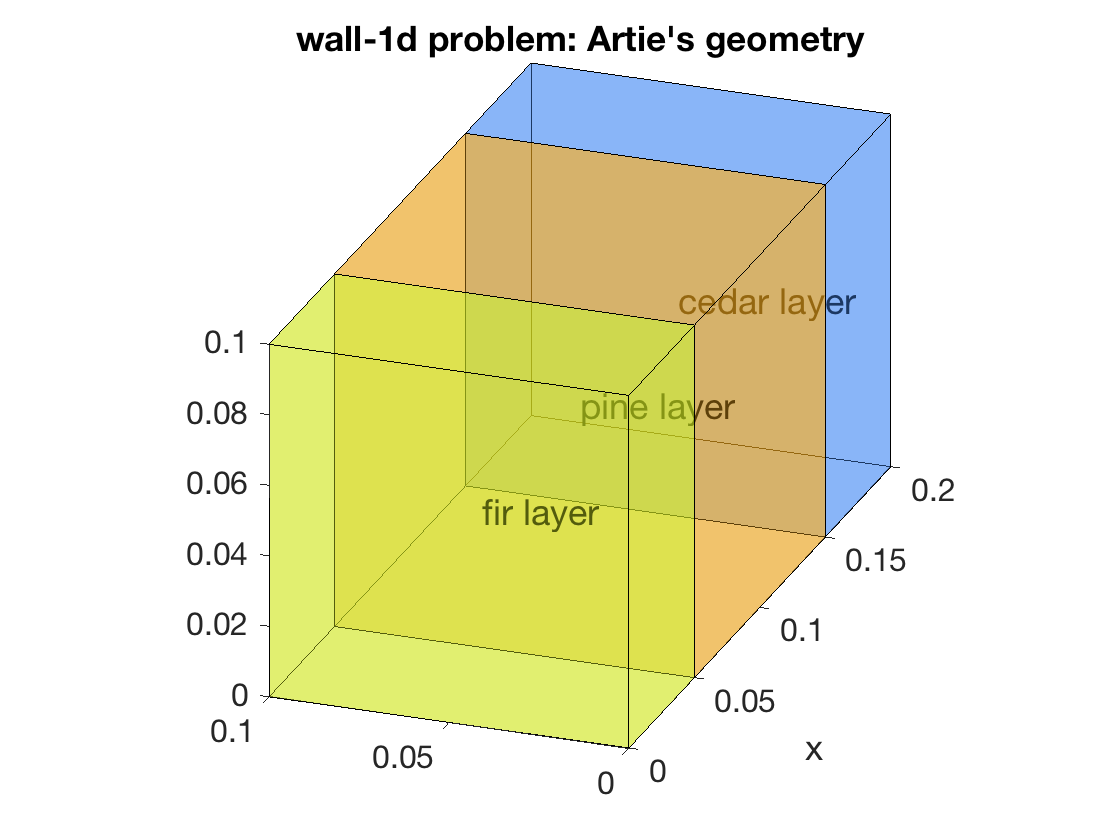}
\end{center}
\caption{Artie's construction of the geometry for problem \texttt{wall-1d}.}
\label{fig:wall-1d_fig2}
\end{figure}

\begin{figure}[h!]
\begin{center}
 \includegraphics[scale=0.3]{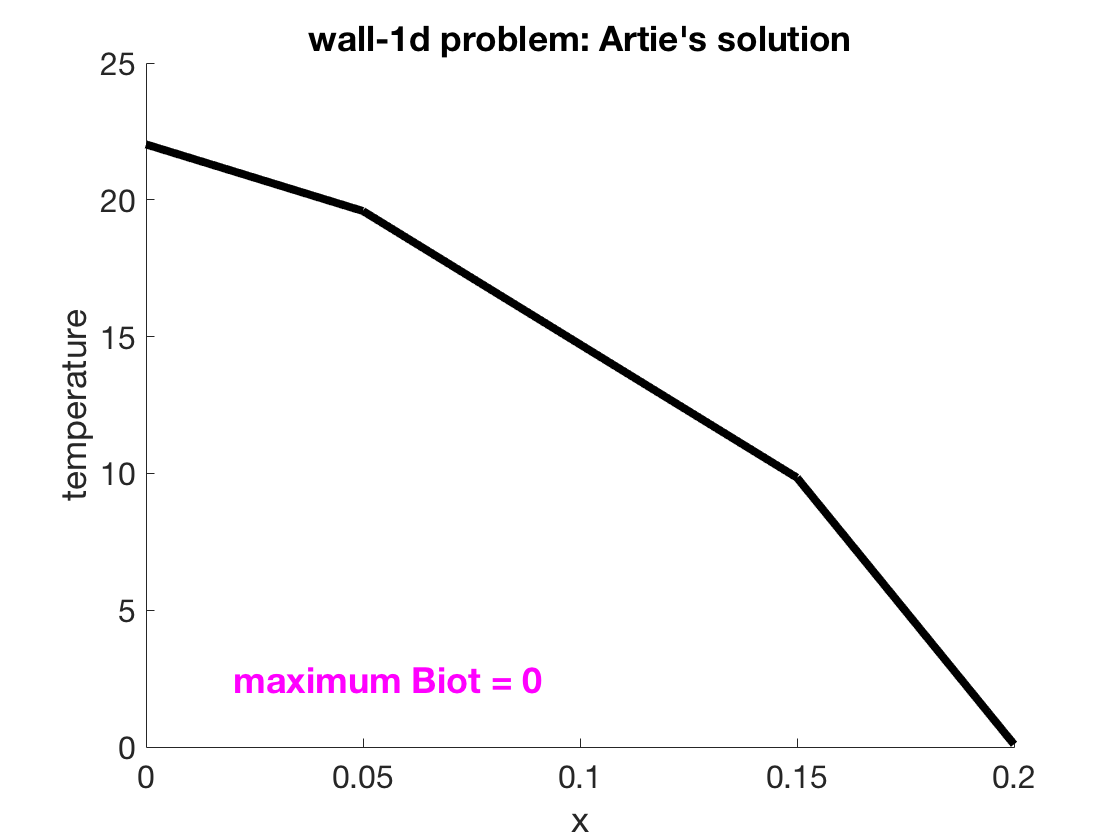}
\end{center}
\caption{Artie's prediction of the temperature field for problem \texttt{wall-1d}.}
\label{fig:wall-1d_fig3}
\end{figure}

Artie's solution is summarized by (i) Figure 1 and Figure 2 referenced in the question statement (but generated by Artie), here respectively Figure \ref{fig:wall-1d_fig1} and Figure \ref{fig:wall-1d_fig2}, and (ii) the plot of the temperature field, here Figure \ref{fig:wall-1d_fig3}. Note this problem has three components corresponding to the three layers of wood; the (children) layers inherit the geometry class of the (parent) wall. The problem is relatively simple both for the students and also for Artie. Note, however, that the actual students and Artie would proceed in different fashions: the student would typically consider a network of convection and conduction resistances in series and apply voltage divider relations; Artie considers the weak formulation with piecewise linear functions within each component.

\subsection{Example 2: \texttt{spoon}}. 

We present the problem statement:

\begin{itemize}
\item[]
A cup contains tea. Air is in contact with the tea; also air surrounds the cup. The tea is maintained at temperature $T_{\rm liq}$; the air is maintained at temperature $T_{\infty}$. The tea - air interface is at $ x = 0 $. A spoon comprises two parts: a head connected to a handle. The head rests in the cup. The head is immersed in the tea; the handle is exposed to the air. Figure 1 is Artie's depiction of the heat transfer paths.

The spoon geometry is approximated as a right cylinder with rectangular cross-section for dimensions $a$ and $b$; the axial coordinate is $x$. The head is of length $L_1$; the handle is of length $L_2$. For our coordinate system, the head extends from $ x = - L_1 $ to $ x = 0 $, and the handle extends from $ x = 0 $ to $ x = L_2 $. Figure 2 is Artie's construction of the spoon geometry.

Let $h_1^{\rm bot}$ denote the heat transfer coefficient from head to tea prescribed over the axial face at $ x = - L_1$; let $h_2^{\rm top}$ denote the heat transfer coefficient from handle to air prescribed over the axial face $x = L_2$; let $h_1^{\rm lat}$ denote the heat transfer coefficient from head to tea prescribed over the head lateral face; let $h_2^{\rm lat}$ denote the heat transfer coefficient from handle to air prescribed over the handle lateral face.

The head has thermal conductivity $k_1$; the handle has thermal conductivity $k_2$; the cup is an insulator.

Plot the temperature as a function of the axial coordinate. You may use the following numerical values: $ a = 0.002$, $b = 0.01 $, $ L_1 =  0.05 $, $L_2 = 0.12$, $h_1^{\rm bot} = 10.0$, $h_2^{\rm top} = 5$, $h_1^{\rm lat} = 10.0$, $h_2^{\rm lat} = 5.0$, $ k_1 = 50 $, $ k_2 = 50$, $T_{\rm liq} = 90$, and $T_\infty = 23$. $\Box$
\end{itemize}
We observe that this problem is in problem class Quasi-1d; the geometry class is right cylinder with rectangular cross-section. We recall that ``plot'' is mandatory within problem class Quasi-1d, and hence the first sentence of the last paragraph is parsed by Artie but no inference is required; the command is primarily for the actual students. We do not request a QoI since in this case the most relevant output --- the temperature of the spoon at $x = L_2$ (as sensed by a finger) --- can be deduced from the field plot.

\begin{figure}[h!]
\begin{center}
 \includegraphics[scale=0.24]{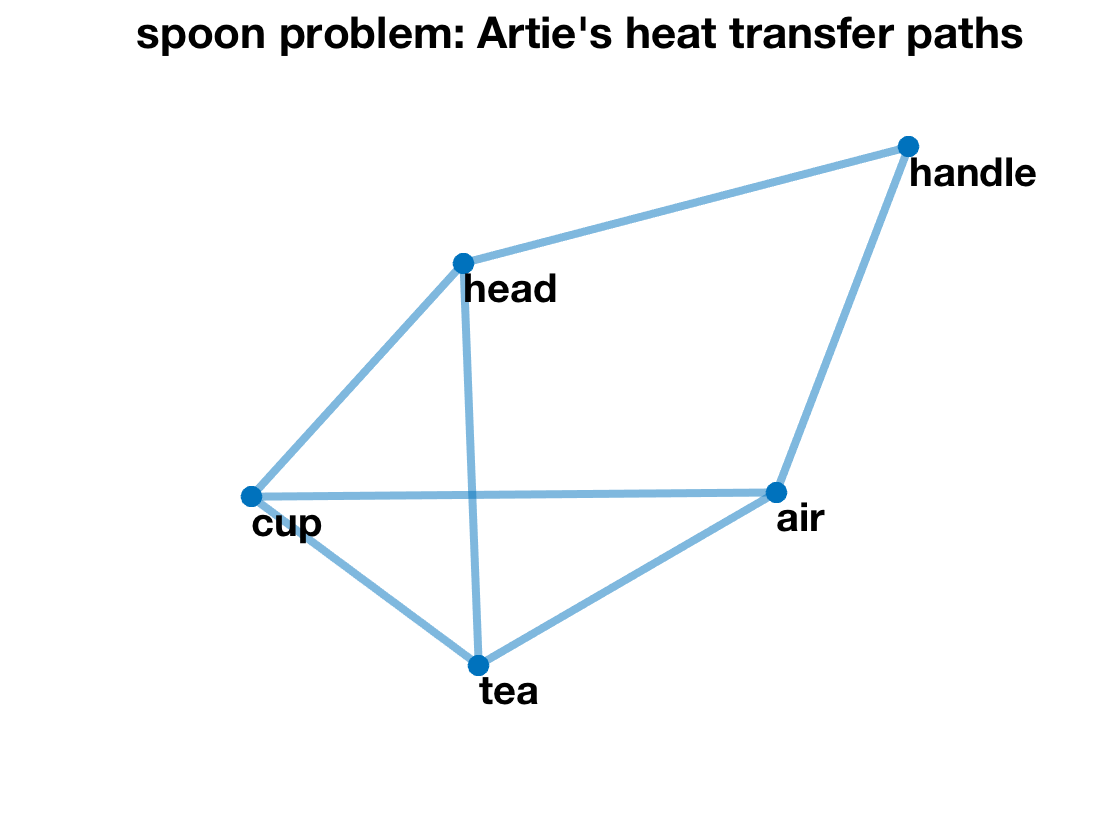}
\end{center}
\caption{Artie's depiction of the heat transfer paths for problem \texttt{spoon}.}
\label{fig:spoon_fig1}
\end{figure}

\begin{figure}[h!]
\begin{center}
 \includegraphics[scale=0.3]{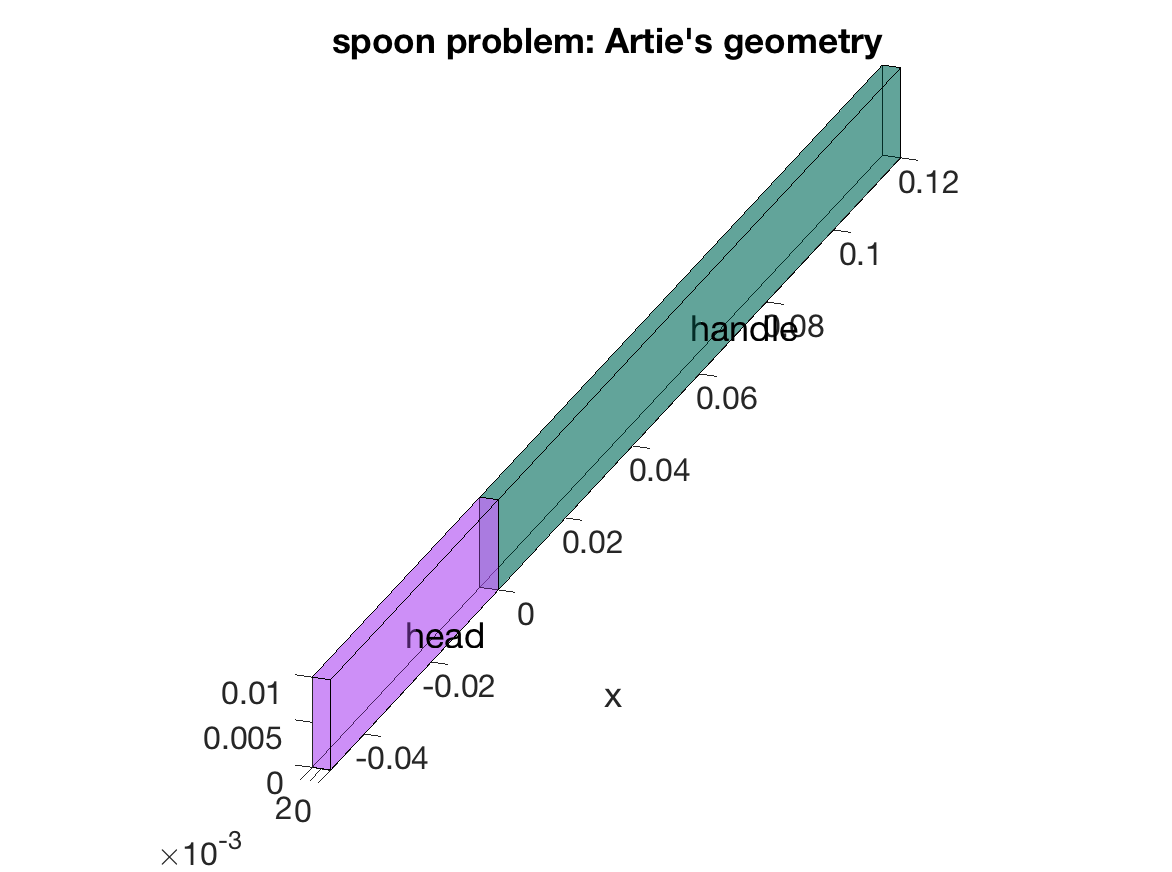}
\end{center}
\caption{Artie's construction of the geometry for problem \texttt{spoon}.}
\label{fig:spoon_fig2}
\end{figure}

\begin{figure}[h!]
\begin{center}
 \includegraphics[scale=0.3]{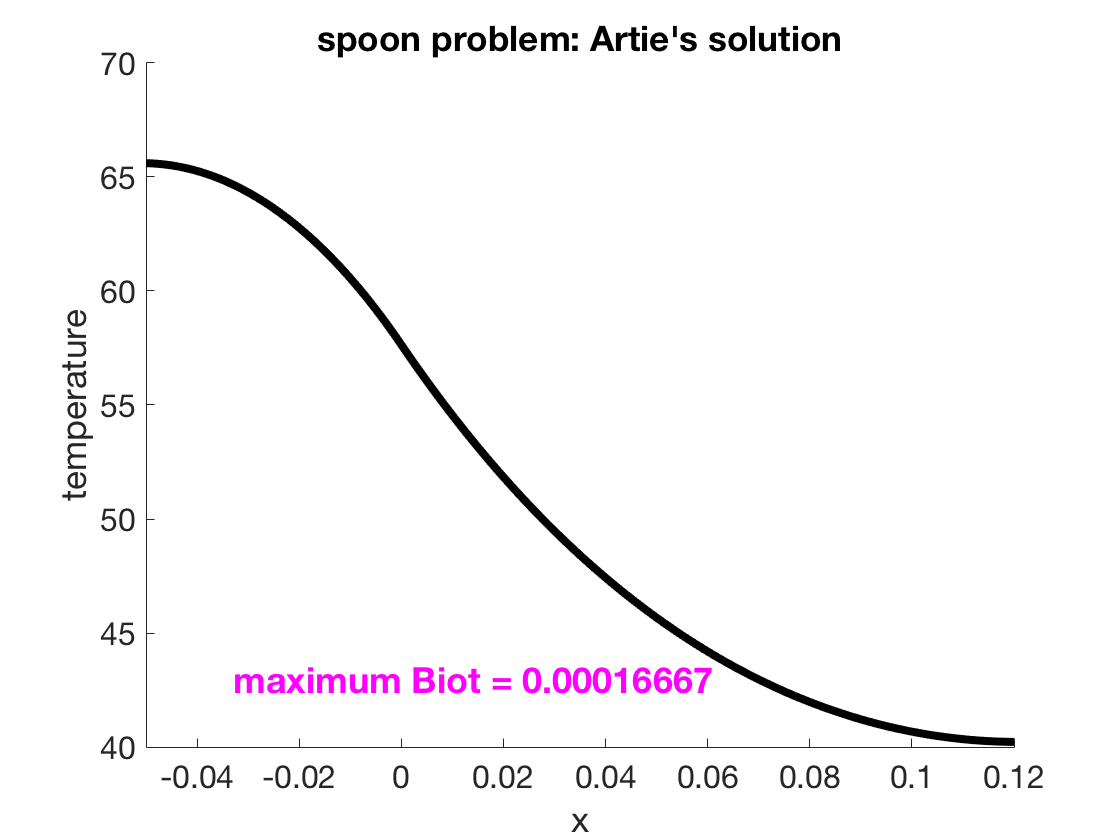}
\end{center}
\caption{Artie's prediction of the temperature field for problem \texttt{spoon}.}
\label{fig:spoon_fig3}
\end{figure}

Artie's solution is summarized by (i) Figure 1 and Figure 2 referenced in the question statement (but generated by Artie), here respectively Figure \ref{fig:spoon_fig1} and Figure \ref{fig:spoon_fig2}, and (ii) the plot of the temperature field, here Figure \ref{fig:spoon_fig3}; note Artie confirms that the Biot number is small, as required by the problem class. This problem has two components corresponding to the head and the handle; the (children) head and handle inherit the geometry class of the (parent) spoon. This problem is less straight-forward: we do not indicate which equation or which limit should be considered. (In fact, most undergraduate students could only solve this problem for the case in which $h_1^{\rm lat} = h_2^{\rm lat}$.)  Again the actual students and Artie would proceed in different fashions: the student would typically invoke the already-derived solution of the 1d fin equation as provided for example in \cite{Lienhard}; Artie considers the weak formulation with appropriate $\sinh$ and $\cosh$ (statically condensed) basis functions, as well as a bubble function, within each component.

\subsection{Example 3: \texttt{wall-3d}}. 

We present the problem statement:

\begin{itemize}
\item[] A wall separates inside air and outside air. The wall consists of four bricks: brick 1, brick 2, brick 3, and brick 4. Brick 1 connects to brick 2; brick 2 also connects to brick 3 and brick 4. Brick 1 is in communication with inside air at temperature $T_{\rm in}$; Brick 2, Brick 3, and Brick 4 are in communication with outside air at temperature $T_{\rm out}$. Figure 1 is Artie's depiction of the heat transfer paths.

The coordinates are denoted $x_1$, $x_2$, and $x_3$; $x_1$ corresponds to distance through the wall. Each brick is a parallelepiped of rectangular cross-section of dimensions $a$ (in $x_2$), $b$ (in $x_3$), $L$ (in $x_1$). The spatial domain of Brick 1 is $0 < x_1 < L, 0 < x_2 < a, 0 < x_3 < b$; the spatial domain of Brick 2 is $L < x_1 < 2L, 0 < x_2 < a, 0 < x_3 < b$; the spatial domain of Brick 3 is $L < x_1 < 2L, 0 < x_2 < a, b < x_3 < 2b$; the spatial domain of Brick 4 is $L < x_1 < 2L, 0 < x_2 < a, - b < x_3 < 0$.  Figure 2 is Artie's construction of the wall geometry.

Each brick has thermal conductivity $k_{\rm b}$.

Brick 1 is exposed to inside air over the face at $x_1 = 0$ through heat transfer coefficient $h_{\rm in}$. Brick 2, Brick 3, and Brick 4 are exposed to outside air over the faces at $x_1 = 2L$ through heat transfer coefficient $h_{\rm out}$. The remainder of the boundary is insulated.

We introduce a nondimensional heat transfer rate $H$ given by $Q/(k_b(T_{\rm in}-T_{\rm out})a)$; here $Q$ denotes the heat transfer rate into Brick 1 over the face at $x_1 = 0$. Develop a lower bound and also an upper bound for $H$. You may use the following parameter values: $T_{\rm in} = 23$, $T_{\rm out} = 0$, $a = 0.1$, $b = 0.1$, $L = 0.05$, $h_{\rm in} = 10$, $h_{\rm out} = 100$, and $k_b = 0.5$. $\Box$
\end{itemize}
We observe that this problem is the problem class Generalized Wall; the geometry class is parallelepiped. The QoI $H$ and associated lower and upper bounds are mandatory with the problem class, so the first two sentences of the last paragraph are parsed by Artie but no inference is required; the QoI request is primarily for the actual students.

\begin{figure}[h!]
\begin{center}
 \includegraphics[scale=0.24]{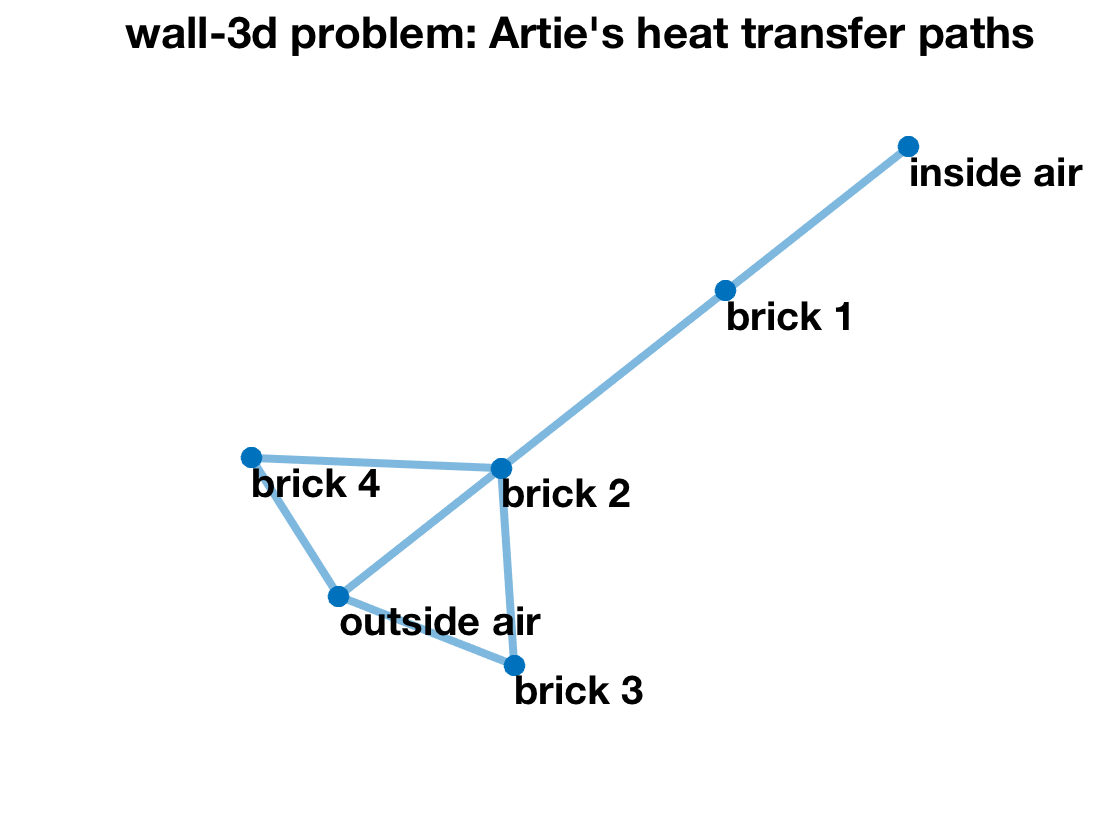}
\end{center}
\caption{Artie's depiction of the heat transfer paths for problem \texttt{wall-3d}.}
\label{fig:wall-3d_fig1}
\end{figure}

\begin{figure}[h!]
\begin{center}
 \includegraphics[scale=0.3]{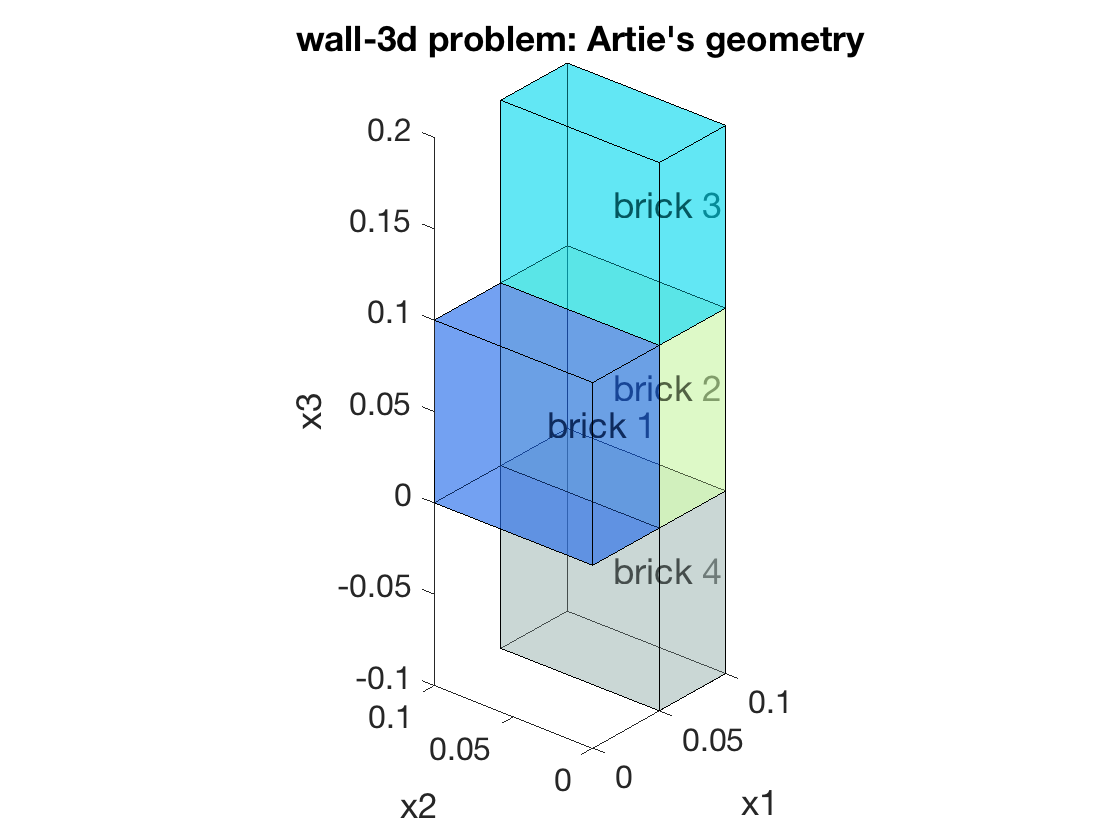}
\end{center}
\caption{Artie's construction of the geometry for problem \texttt{wall-3d}.}
\label{fig:wall-3d_fig2}
\end{figure}

\begin{figure}[h!]
\vspace{.2in}
\begin{mdframed}
\begin{center}
 \includegraphics[scale=0.7]{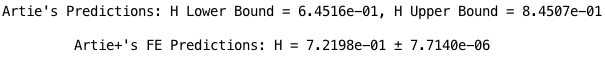}
\end{center}
\end{mdframed}
\vspace{-.1in}
\caption{Artie's lower and upper bound for $H$ for problem \texttt{wall-3d}; Artie+'s FE prediction for $H$ and associated {\it a posteriori} error estimator.}
\label{fig:wall-3d_fig4}
\end{figure}

\begin{figure}[h!]
\begin{center}
 \includegraphics[scale=0.3]{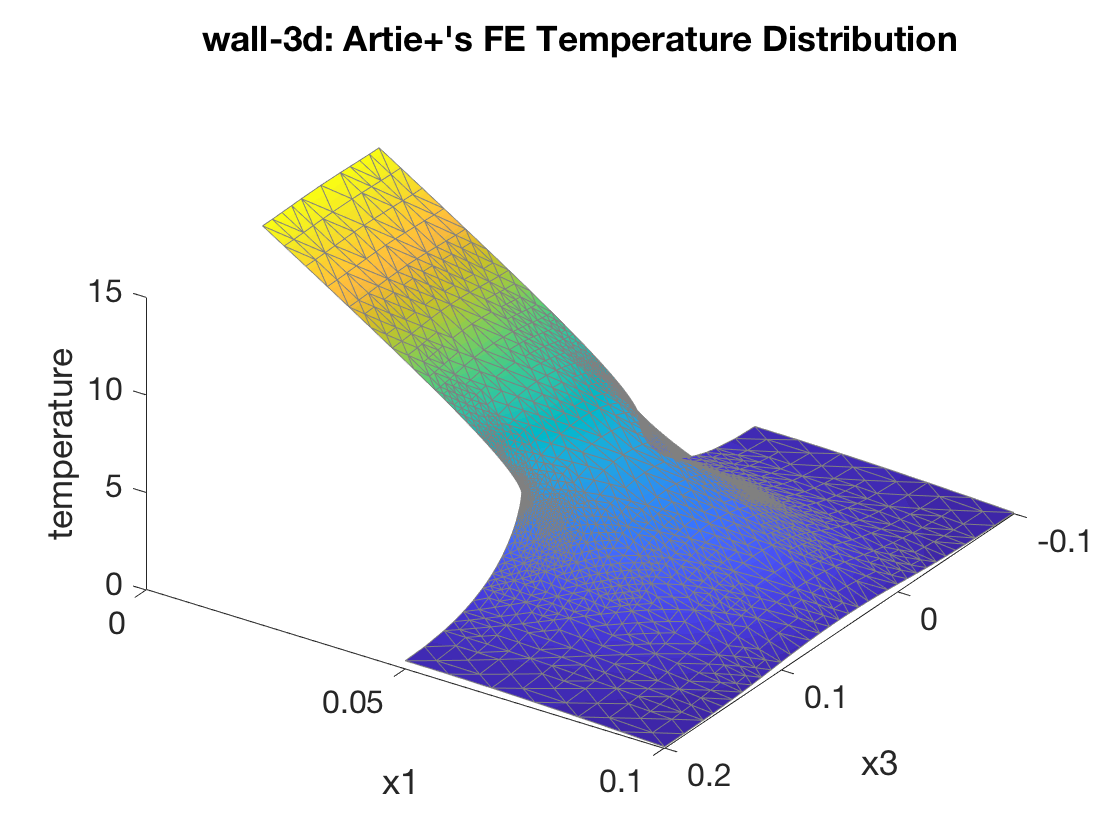}
\end{center}
\caption{Artie+'s FE prediction of the temperature field for problem \texttt{wall-3d} over any $x_2 = \text{ const}$ plane.}
\label{fig:wall-3d_fig3}
\end{figure}

Artie+'s solution is summarized by (i) Figure 1 and Figure 2 referenced in the question statement (but generated by Artie), here respectively Figure \ref{fig:wall-3d_fig1} and Figure \ref{fig:wall-3d_fig2}, (ii) the display of the lower and upper bounds for $H$, here provided in Figure \ref{fig:wall-3d_fig4}, and (iii) the plot of the FE temperature field, here Figure \ref{fig:wall-3d_fig3}. As regards the latter, the geometry is 3d but insulated in $x_2$ and hence the temperature field is in fact 2d and thus amenable to our (currently) two-dimensional finite element treatment. This problem has four components corresponding to Brick $n$, $1 \le n \le 4$; Brick $n$, $1 \le n \le 4$, are instantiations of the archetype brick and inherit from the archetype the geometry class and also conductivity (any conductivities provided explicitly in the problem statement for a particular Brick $n'$ would take precedence).\footnote{Note the \texttt{wall-3d} problem statement is somewhat inconsistent as regards capitalization of proper bricks; Artie is not confused.} This problem is not too difficult for students {\em if} we also include a hint --- as we would typically do in practice --- for the bound construction: for the lower bound, insulator cuts are inserted on $\bar \Omega \cap \{x_3 = 0\}$ and $\bar \Omega \cap \{x_3 = b\}$, whereas for the upper bound, a superconductor cut is inserted on $\bar \Omega \cap \{x_1 = L\}$; in both cases, the resulting heat flow is 1d and may be treated by resistances in series and parallel. (Note $\{x_1 = \text{ const}\}$ is shorthand for $\{ x \in \RR^3 \,| \, x_1 = \text{ const}\}$.) We need not provide these hints to Artie; Artie can deduce the bound construction from the connectivity information.

We elaborate upon Artie's approach and in particular the variational formulation. In this analysis we assume that $h_{\rm in} > 0$ and $h_{\rm out} > 0$, and also, as always, $k_{\rm b} > 0$. We also introduce several standard function spaces: $L^2(\Omega)$ is the space of functions which are (Lebesgue) square integrable over $\Omega$; $H^1(\Omega) \equiv \{ v \in L^2(\Omega)\,|\, |\nabla v| \in L^2({\Omega})\}$. For our problem (and more generally problem class Generalized Wall) we can evaluate $H$ as
\begin{align}
H = \frac{1}{k_{\rm b}(T_{\rm in} - T_{\rm out})a}\int_{\partial \Omega_0} h_{\rm in} (T_{\rm in} - T)\,dA \;,  \label{H1}
\end{align}
where $T$ refers to the solution of our PDE for the prescribed boundary conditions; note $\partial \Omega_0 \equiv \{ x \in \partial \Omega \, | \, x_1 = 0\}$ and, for future reference, $\partial \Omega_{2L} \equiv \{ x \in \partial \Omega \, | \, x_1 = 2L\}$. We can also express $H$ as a minimum:
\begin{align}
H = C_1(\min_{w \in X} 2 J(w) + C_2)\; , \label{H2}
\end{align}
where $C_1$ and $C_2$ are positive constants,
\begin{align}
C_1 = \frac{1}{k_{\rm b}(T_{\rm in} - T_{\rm out})^2 a}\; , \label{C1}
\end{align}
\begin{align}
C_2 =  \int_{\partial \Omega_0} h_{\rm in} T_{\rm in}^2 + \int_{\partial \Omega_{2L}} h_{\rm out} T_{\rm out}^2 \; , \label{C2}
\end{align}
$X = H^1(\Omega)$, and $J: X \times X \rightarrow \RR$ is the functional
\begin{align}
J(w) = \frac{1}{2} \calA(w,w) - \calF(w) \; ; \label{J}
\end{align}
here $\calA$ and $\calF$ are respectively the (symmetric, coercive) bilinear and linear forms associated to our PDE,
\begin{align}
\calA(w,v) \equiv \int_{\partial \Omega_0} h_{\rm in} w v\, dA + \int_{\partial \Omega_{2L}} h_{\rm out} w v dA + \int_\Omega k_{\rm b} \nabla w \cdot \nabla v\,dV  \; , \label{calA}
\end{align}
and
\begin{align}
\calF(v) \equiv \int_{\partial \Omega_0} h_{\rm in} T_{\rm in} v\, dA + \int_{\partial \Omega_{2L}} h_{\rm out} T_{\rm out} v\, dA\; . \label{calF}
\end{align}
We note that $T \in X$ satisfies the weak form: $\calA(T,v) = \calF(v), \forall v \in X$.

Artie now proceeds to identify the $x_1$ locations of the brick (component) $x_1$ faces in the system, $\{x_{1}^j\}_{j = 1,\ldots,m_1}$, and associated domain slices $\calS^j \equiv \text{ (the interior of) } \bar \Omega \cap \{x_1 = x_1^j\}$; for \texttt{wall-3d},  $m_1 = 3$ with $x_1^j = (j-1)L, 1 \le j \le m_1$, and $\calS^1 \equiv \{x_1 = 0, 0 < x_2 < a, 0 < x_3 < b\}$, $\calS^2 \equiv \{x_1 = L, 0 < x_2 < a, -b < x_3 < 2b\}$, $\calS^3 \equiv \{x_1 = 2L, 0 < x_2 < a, -b < x_3 < 2b\}$. We also identify a minimal set of parallelepipeds, $\{\calP^j\}_{j = 1,\ldots,m_{\rm ppd}}$, such that $\bar \Omega = \cup_{j = 1}^{m_{\rm ppd}} \bar \calP^j$: each parallelepiped is a (sub)assembly of bricks with connections only on brick (component) $x_1$ faces; for \texttt{wall-3d}, $m_{\rm ppd} = 3$ with $\calP^1 \equiv \{0 < x_1 < 2L, 0 < x_2 < a, 0 < x_3 < b\}$, $\calP^2 \equiv \{L < x_1 < 2L, 0 < x_2 < a, b < x_3 < 2b\}$, $\calP^3 \equiv \{L < x_1 < 2L, 0 < x_2 < a, -b < x_3 < 0\}$. We then introduce two additional spaces: 
\begin{align}
X^{\rm UB}= \{ v \in H^1(\Omega)\,| \, v|_{\calS^j} = \text{ Const}^j, j = 1,\ldots,m_1\} \; , \label{XUB}
\end{align}
and
\begin{align}
X^{\rm LB} = \{ v \in L^2(\Omega)  \, | \, v|_{\calP^j}  \in H^1(\calP^j), j = 1,\ldots,m_{\rm ppd}\} \; . \label{XLB}
\end{align}
We make two remarks on the upper-bound space, \eqref{XUB}. First, for any given $j$, the constant Const$^j \in \RR$ that appears in the constraint is {\em part} of the solution, hence unknown {\it a priori} --- the constraint imposes uniformity over $\calS^j$, {\em not} a prescribed (Dirichlet) value. Second, if any $\calS^j$ is disconnected --- not the case for \texttt{wall-3d} --- we can improve the upper bound space: we permit not a single constant Const$^j$ associated to the entire slice $\calS^j$ but rather a different constant for each {\em connected part} of slice $\calS^j$. We also make a remark on the lower-bound space, \eqref{XLB}: for any parallelepiped $\calP$ with no Robin  faces we must additionally require that the field satisfy a zero-mean condition over $\calP$. We can then construct the quantities
\begin{align}
H^{\rm UB} = C_1(\min_{w \in X^{\rm UB}}2J(w) + C_2) \; , \label{HUB}
\end{align}
and
\begin{align}
H^{\rm LB} = C_1(\min_{w \in X^{\rm LB}} \sum_{j = 1}^{m_{\rm ppd}} 2 J|_{\calP^{j}}(w) + C_2) \; , \label{HLB}
\end{align}
which, from simple variational arguments, satisfy $H^{\rm LB} \le H \le H^{\rm UB}$. Note that $J|_{\calP^j}$ refers to $J$ restricted to $\calP^j$.

The implementation of the bounds is very simple, and thus within Artie's reach even without any knowledge of variational methods. We discuss the upper bound; the lower bound admits even simpler treatment. The exact minimizer associated with \eqref{HUB} is in fact piecewise-linear in $x_1$ (and independent of $x_2$ and $x_3$); Artie can thus deduce the discrete equations directly from the weak form over $X^{\rm UB}$ applied to basis functions which are piecewise linear over each component. The only subtlety, easily addressed, is the constraint in $X^{\rm UB}$: for each $j \in \{1,\ldots,m_1\}$, the global degrees of freedom associated to all ports on  (component) $x_1 = x_1^j$ faces are simply coalesced into a single degree of freedom --- a revised local to global mapping from component and local port to unique global degree of freedom --- prior to direct stiffness assembly over components.

\section{Future Work}

There are a number of short-term targets of opportunity: improvement of Artie's parsing skills in particular for compound predicates; generation of more examples including different boundary conditions and various QoI; development of new geometry classes such as (for problem class Quasi-1d) axisymmetric and spherical configurations, and more generally small-curvature filaments; extension to other fields of heat transfer (unsteady conduction, radiation), to linear elasticity and structural mechanics, and to fluid dynamics and heat convection.

We also envision more challenging tasks. First, we would like to extend Artie to consider {\em problem perturbation} as discussed in Pedagogical Background, Section \ref{subsec:PB}: given a problem statement, find a ``closest'' problem class and associated parameters; the latter is best conducted in concert with the development of {\em image processing} techniques to deduce connections and geometry. Important considerations include the choice of norm and the scoring algorithm; QoI bounds can in certain cases play an important role. Second, we would like to initiate a first simple form of learning, {\em commonsense learning}: Artie can pose questions (``Is tea a solid?") the answer to which can then be incorporated into the commonsense database for use in future problems. Third, we would like to incorporate component-based {\em model-order reduction} techniques \cite{Enu} into Artie+'s PDE numerical solution capabilities in particular to reduce response time for truly three-dimensional problems.

Finally, in the much longer term, we would need to develop an uber-Artie who could develop Artie software for a particular discipline from associated textbook material: {\em textbook learning}. We would continue to rely on the PDE Template but now at a higher level: the PDE Template would serve to guide Artie's knowledge acquisition and subsequent development of discipline-specific parsers. In the even longer term, Artie could incorporate archival material for example related to material properties, heat transfer coefficient correlations, and other empirical information.

\section{Acknowledgments}

I would like to thank Dr Phuong Huynh of Akselos, SA, for fruitful exchanges about Project Artie and in particular for his development and generous distribution of the \texttt{GNLRequest} function for direct client access to the Google Natural Language Processor through MATLAB;  the latter greatly facilites syntax analysis tightly integrated with symbolic and numeric computation. I would like to thank Professor Masayuki Yano of University of Toronto for his pedagogical and mathematical insights into Project Artie and furthermore for the development and generous distribution of his \texttt{fem2d}  templates for  adaptive finite element solution of partial differential equations. I would also like to thank Professor Ed Greitzer of MIT for his comments on Artie in education, and Professor Reza Malek-Madani of the US Naval Academy for his encouragement to pursue Project Artie. Project Artie is supported by MIT through a Ford Professorship.


\footnotesize
\bibliographystyle{plain}	%
\bibliography{bib_folder/Artie}

\newpage
\appendix
\section{Problem Statements: Inputs to Artie+}

\includepdf[pages=1,pagecommand=\subsection{Example 1: \texttt{wall-1d}}, offset = 0 -3cm]{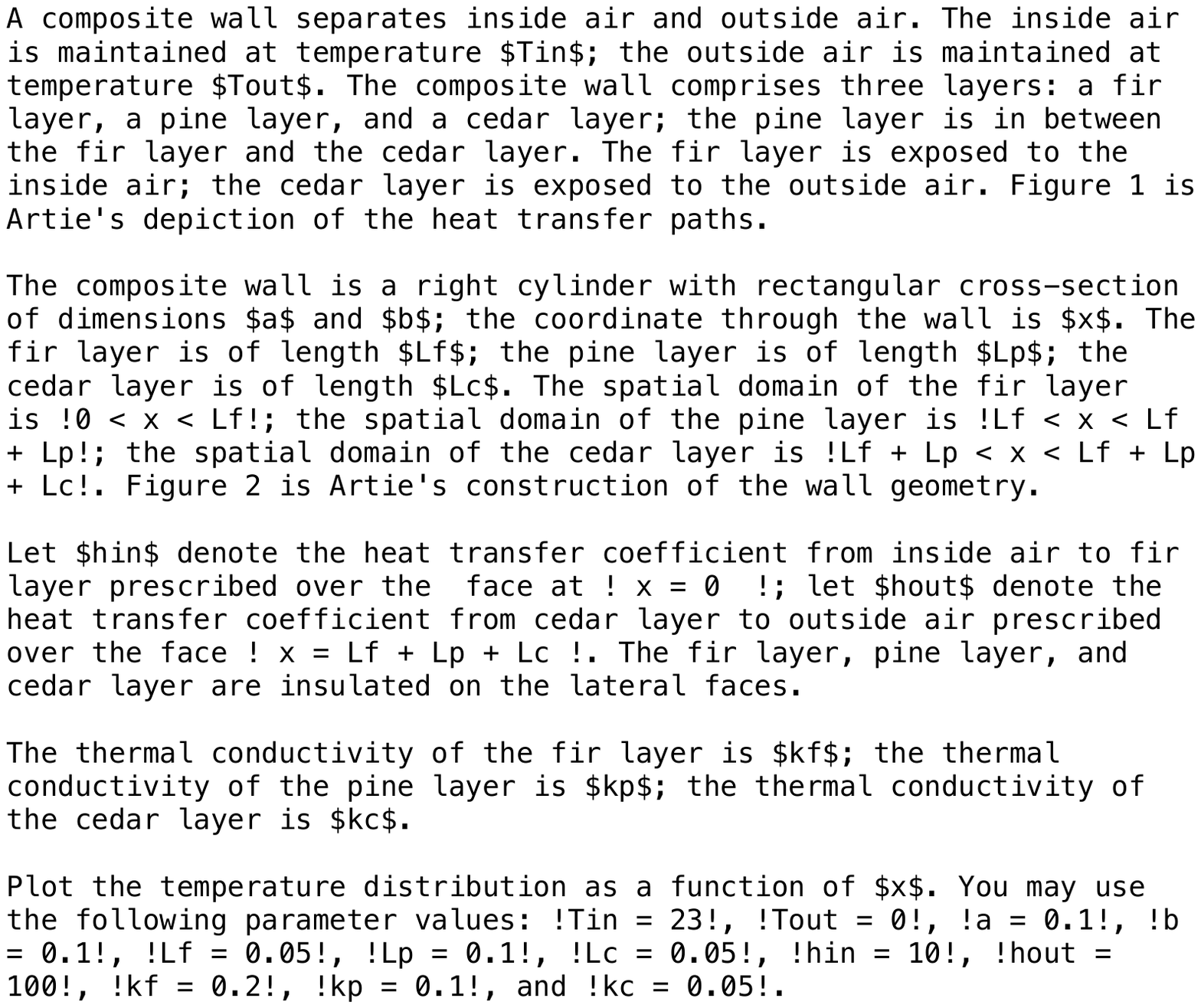}
\includepdf[pages=1,pagecommand=\subsection{Example 2: \texttt{spoon}}, offset = 0 -3cm]{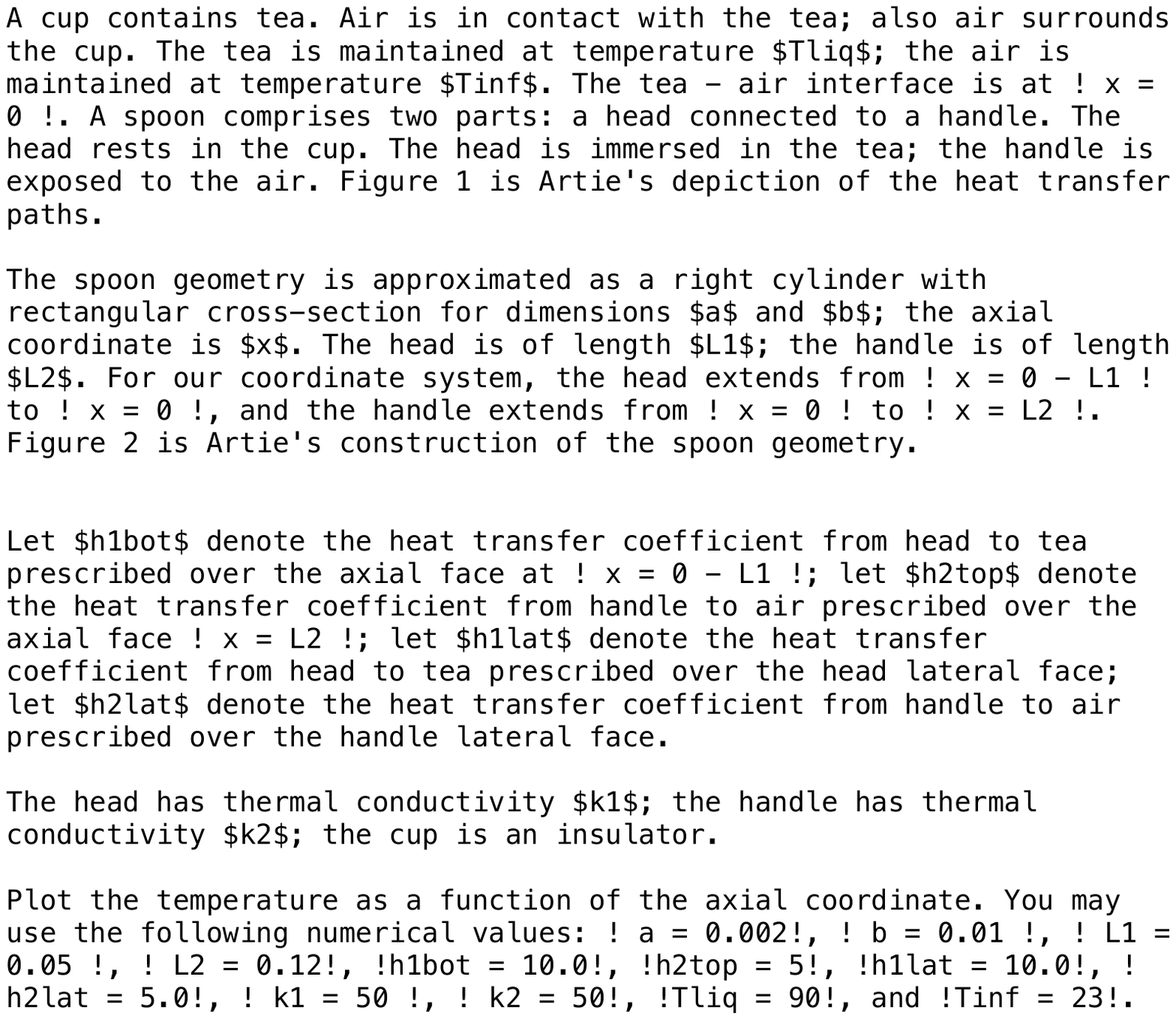}
\includepdf[pages=1,pagecommand=\subsection{Example 3: \texttt{wall-3d}}, offset = 0 -3cm]{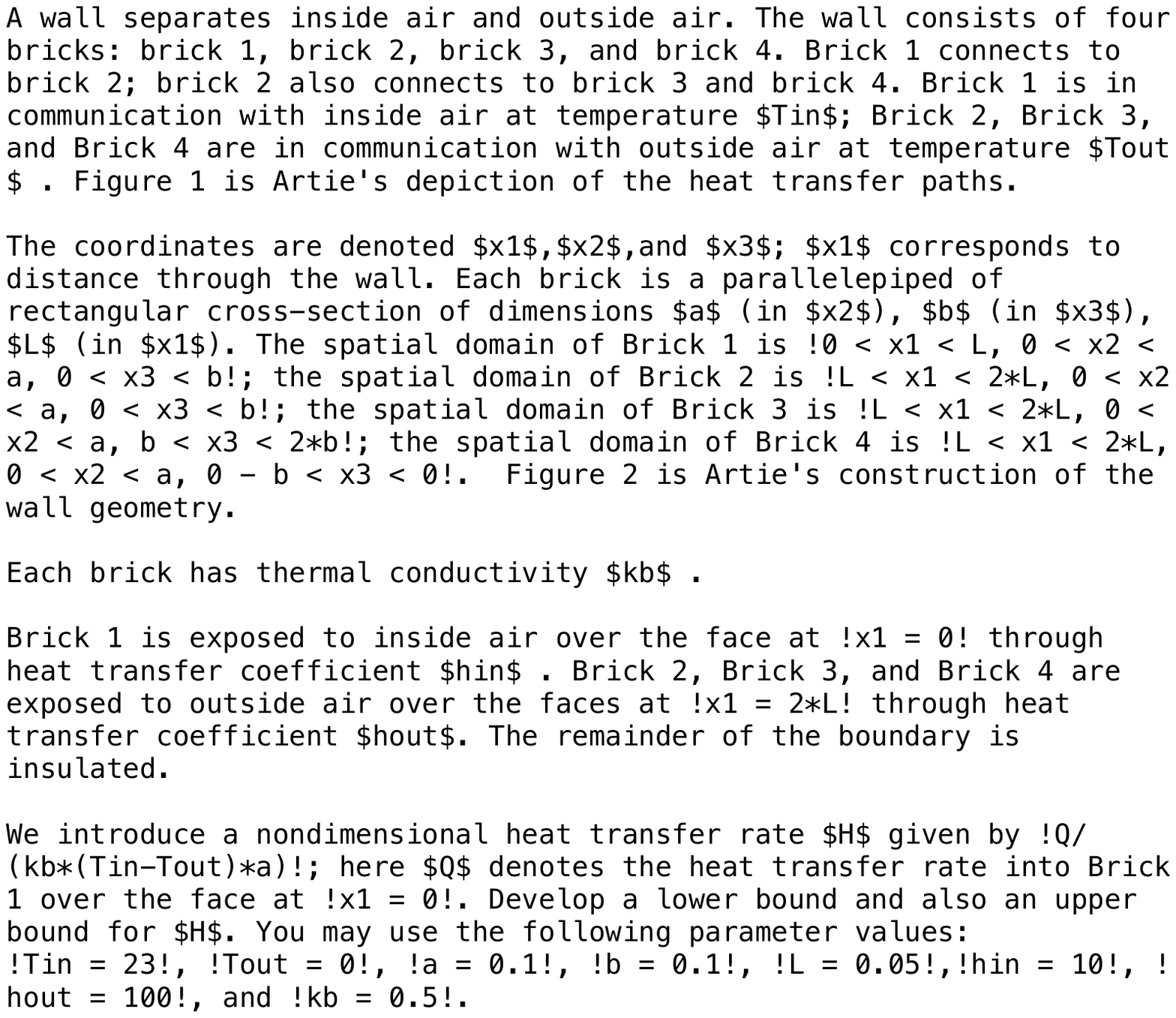}


%

\end{document}